\documentclass[11pt]{amsart}
\usepackage{color}
\usepackage{amssymb}
\usepackage{amsmath}

\newcommand{\bbH}{{\mathbb H}}
\newcommand{\bbZ}{{\mathbb Z}}

\newcommand{\bbP}{{\mathbb P}}

\newcommand{\Oh}{{\mathcal O}}

\DeclareMathOperator{\HH}{H}

\DeclareMathOperator{\Pic}{Pic}

\DeclareMathOperator{\Ext}{Ext}

\DeclareMathOperator{\rank}{rank}

\DeclareMathOperator{\pd}{pd}

\newcommand{\ol}[1]{\overline{#1}}
\newcommand{\onto}{\twoheadrightarrow}

\newcommand{\by}[1]{\xrightarrow{#1}}

\newcommand{\tensor}{\otimes}
\newcommand{\isom}{\cong}

\newcommand{\sE}{{\mathcal E}}

\newcommand{\sExt}{\mbox{${\sE}{xt}$}}

\theoremstyle{plain}
\newtheorem{lemma}{Lemma}
\newtheorem{prop}{Proposition}

\newtheorem{cor}{Corollary}
\newtheorem*{question}{Questions}
\newtheorem*{question1}{Question}

\newtheorem*{ACM}{ACM Degree Conjecture}
\theoremstyle{definition}

\DeclareMathOperator{\depth}{depth}
\long\def\comment#1{}
\begin{document}

\title{On codimension two subvarieties in hypersurfaces}
\author{N. Mohan Kumar}
\address{Department of Mathematics, Washington University in
St. Louis, St. Louis, Missouri, 63130}
\email{kumar@wustl.edu}
\urladdr{http://www.math.wustl.edu/$\sim$kumar}
\author{A.~P.~Rao}
\address{Department of Mathematics, University of Missouri-St. Louis,
St. Louis, Missouri 63121}
\email{rao@arch.umsl.edu}
\author{G. V. Ravindra}
\address{Department of Mathematics, Indian Institute of Science,
Bangalore 560012, India}
\email{ravindra@math.iisc.ernet.in}
\thanks{We thank the referee for pointing out some relevant references.}
\begin{abstract}
We show that for a smooth hypersurface $X\subset \bbP^n$ of degree at
least $2$, there exist arithmetically Cohen-Macaulay (ACM) codimension
two subvarieties $Y\subset X$ which are not an intersection $X\cap{S}$
for a codimension two subvariety $S\subset\bbP^n$. We also show there
exist $Y\subset X$ as above for which the normal bundle sequence for
the inclusion $Y\subset X\subset\bbP^n$ does not split.
\end{abstract}
\date{\today}
\maketitle

\centerline{\it Dedicated to Spencer Bloch}

\section{Introduction}
In this note, we revisit some questions  of Griffiths
and Harris from 1985 \cite{GH}:

\begin{question}[Griffiths and Harris]\label{ghconj}
Let $X\subset\bbP^4$ be a general hypersurface of degree $d\geq 6$ and
$C\subset X$ be a curve. 
\begin{enumerate}
\item \label{ghconj3}Is the degree of $C$  a multiple of $d$?
\item \label{ghconj1}Is $C=X\cap S$ for some surface $S\subset\bbP^4$?
\end{enumerate}\end{question}

The motivation for these questions comes from trying to extend the
Noether-Lefschetz theorem for surfaces to threefolds. Recall that the
Noether-Lefschetz theorem states that if $X$ is a very general surface
of degree $d\geq 4$ in $\bbP^3$, then $\Pic(X)=\bbZ$, and hence every
curve $C$ on $X$ is the complete intersection of $X$ and another
surface $S$.

C.~Voisin very soon \cite{Vo} proved that the second question had a
negative answer by constructing counter-examples on any smooth
hypersurface of degree at least 2.  She also considered a third
question:

\begin{question1} With the same terminology and when $C$ is smooth:
\begin{enumerate}
\setcounter{enumi}{2}
\item \label{ghconj2} Does the exact sequence of normal bundles
  associated to the inclusions $C\subset X\subset\bbP^4$:
$$ 0 \to N_{C/X} \to N_{C/\bbP^4} \to \Oh_C(d) \to 0 $$
split?
\end{enumerate}
\end{question1}

Her counter-examples provided a negative answer to this question as
well.  The first question, the Degree Conjecture of Griffiths-Harris,
is still open.  \comment{\textcolor[rgb]{0.98,0.00,0.00}}Strong
evidence for this conjecture was provided by some elementary but
ingenious examples of Koll{\'a}r (\cite{BCC},Trento examples). In
particular he shows that if $\gcd(d,6)=1$ and $d\geq 4$ and $X$ is a
very general hypersurface of degree $d^2$ in $\bbP^4$, then every
curve on $X$ has degree a multiple of $d$. In the same vein, Van
Geemen shows that if $d>1$ is an odd number and $X$ is a very general
hypersurface of degree $54d$, then every curve on $X$ has degree a
multiple of $3d$.

The main result of this note is the existence of a large class of
counterexamples which subsumes Voisin's counterexamples and places
them in the context of arithmetically Cohen-Macaulay (ACM) vector
bundles on $X$. It is well known that ACM bundles which are not sums
of line bundles can be found on any hypersurface of degree at least 2
\cite{BGS}, and for such a bundle, say of rank $r$, on $X$, ACM
subvarieties of codimension two can be created on $X$ by considering
the dependency locus of $r-1$ general sections. These subvarieties
fail to satisfy Questions \ref{ghconj1} and \ref{ghconj2}. We will be
working on hypersurfaces in $\bbP^n$ for any $n \geq 4$ and our
constructions of ACM subvarieties may not give smooth ones.  Hence in
Question \ref{ghconj2}, we will consider the splitting of the conormal
sheaf sequence instead.

\section{Main results}
Let $X\subset\bbP^n$ be a smooth hypersurface of degree $d\geq 2$ and
let $Y\subset X$ be a codimension $2$ subscheme. Recall that $Y$ is
said to be \comment{\textcolor[rgb]{0.98,0.00,0.00}}an arithmetically
Cohen-Macaulay (ACM)
\comment{\textcolor[rgb]{0.98,0.00,0.00}}subscheme of $X$ if
$\HH^i(X,I_{Y/X}(\nu))=0$ for $0 < i\leq \dim{Y}$ and for any
$\nu\in\bbZ$. Similarly, a vector bundle $E$ on $X$ is said to be ACM
if $\HH^i(X,E(\nu))=0$ for $i\neq 0, \,\dim{X}$ and for any
$\nu\in\bbZ$.

Given a coherent sheaf $\mathcal{F}$ on $X$, let
$s_i\in\HH^0(\mathcal{F}(m_i))$ for $1\leq i \leq k$ be generators for
the $\oplus_{\nu\in\bbZ}\HH^0(\Oh_X(\nu))$-graded module
$\oplus_{\nu\in\bbZ}\HH^0(\mathcal{F}(\nu))$. These sections give a
surjection of sheaves $\oplus_{i=1}^k\Oh_X(-m_i) \onto \mathcal{F}$
which induces a surjection of global section
$\oplus_{i=1}^k\HH^0(\Oh_X(\nu-m_i)) \onto \HH^0(\mathcal{F}(\nu))$
for any $\nu\in\bbZ$.

Applying this to the ideal sheaf $I_{Y/X}$ of an ACM subscheme of
codimension $2$ in $X$, we obtain the short exact sequence
$$ 0 \to G \to \oplus_{i=1}^k\Oh_X(-m_i) \to I_{Y/X} \to 0,$$ where
$G$ is some ACM sheaf on $X$ of rank $k-1$.
\comment{\textcolor[rgb]{0.98,0.00,0.00}}Since $Y$ is ACM as a
subscheme of $X$, it is also ACM as a subscheme of $\bbP^n$. In
particular, $Y$ is locally Cohen-Macaulay. Hence $G$ is a vector
bundle by \comment{\textcolor[rgb]{0.98,0.00,0.00}} the
Auslander-Buchsbaum Theorem (see \cite{Mat} page 155). We will loosely
say that $G$ is associated to $Y$.

Conversely, the following Bertini type theorem which goes back to
arguments of Kleiman in \cite{Kleiman} (see also \cite{Ban}) shows
that given an ACM bundle $G$ on $X$, we can use $G$ to construct ACM
subvarieties $Y$ of codimension $2$ in $X$:

\begin{prop}\label{Kleiman} (Kleiman). 
Given a bundle $G$ of rank $k-1$ on $X$, a general map $G \to
\oplus_{i=1}^k\Oh_X(m_i)$ for sufficiently large $m_i$ will determine
the ideal sheaf (up to twist) of a subvariety $Y$ of codimension $2$
in $X$ with a resolution of sheaves:
$$ 0 \to G \to \oplus_{i=1}^k\Oh_X(m_i) \to I_{Y/X}(m) \to 0.$$
\end{prop}

Since the conclusion of Question \ref{ghconj1} implies that of
Question \ref{ghconj2}, we will look at just Question \ref{ghconj2},
in the conormal sheaf version.

\medskip

Let $X$ be a hypersurface of degree $d$ in $\bbP^n$ defined by the
equation $f=0$.  Let $X_2$ be the thickening of $X$ defined by $f^2=0$
in $\bbP^n$.  Given a subvariety $Y$ of codimension $2$ in $X$, let
$I_{Y/\bbP}$ (resp. $I_{Y/X}$) denote the ideal sheaf of
$Y\subset\bbP^n$ (resp. $Y\subset X$).  The conormal sheaf sequence is
 \begin{equation}\label{nbs}
0\to \Oh_Y(-d) \to I_{Y/\bbP}/I_{Y/\bbP}^2 \to I_{Y/X}/I_{Y/X}^2 \to 0.
\end{equation}

\begin{lemma}
For the inclusion $Y\subset X\subset \bbP^n$, if the sequence of
conormal sheaves (\ref{nbs}) splits, then there exists a subscheme
$Y_2\subset X_2$ containing $Y$ such that
$$   I_{Y_2/X_2}(-d) \stackrel{f}{\to} I_{Y_2/X_2}  \to I_{Y/X} \to 0$$
is exact. Furthermore, $fI_{Y_2/X_2}(-d)=I_{Y/X}(-d)$.
\end{lemma}

\begin{proof}
Suppose sequence (\ref{nbs}) splits: then we have a surjection
$$ I_{Y/\bbP}\onto I_{Y/\bbP}/I_{Y/\bbP}^2  \onto \Oh_Y(-d)$$ where
the first map is the natural quotient map and the second is the
splitting map for the sequence. The kernel of this composition
defines a scheme $Y_2$ in $\bbP^n$. Since this kernel $I_{Y_2/\bbP}$  
contains $I_{Y/\bbP}^2$ and hence $f^2$, it is clear that
$Y \subset Y_2 \subset X_2$.

The splitting of (\ref{nbs}) also means that $f \in I_{Y/\bbP}(d)$
maps to $1 \in \Oh_Y$.  We get the commutative diagram:

\[
\begin{array}{ccccccccc}
& & & & & & 0 & \\
& & & & & & \uparrow & & \\
0 &\to & I_{Y_2/\bbP}& \to &  I_{Y/\bbP} & \to & \Oh_Y(-d) & \to & 0 \\
& &  \uparrow{\scriptstyle{f^2}} & &  \uparrow{\scriptstyle{f}} & &
\uparrow & &\\ 
0 & \to & \Oh_\bbP(-2d) & \stackrel{f}{\to}
 & \Oh_\bbP(-d) & \to & \Oh_X(-d) &\to & 0\\
& & \uparrow & & \uparrow & & & &\\
& & 0 & & 0 & & & &\\
\end{array}
\]

This induces 
$$ 0 \to I_{Y/X}(-d) \to I_{Y_2/X_2} \to I_{Y/X} \to 0.$$ 
In particular, note that $I_{Y/X}(-d)$ is the image of the multiplication
map $f: I_{Y_2/X_2}(-d) \to I_{Y_2/X_2}$.
\end{proof} 

\medskip
Now assume that $Y$ is an ACM subvariety on $X$ of codimension $2$.
The ideal sheaf of $Y$ in $X$ has a resolution
$$ 0 \to G \to \oplus_{i=1}^k\Oh_X(-m_i) \to I_{Y/X} \to 0,$$ for some 
ACM bundle $G$ on $X$ associated to $Y$.

\begin{lemma}\label{bundle extends} 
Suppose the conditions of the previous lemma hold, and in addition $Y$
is an ACM subvariety. Then there is an extension of the ACM bundle $G$
(associated to $Y$) on $X$ to a bundle $\mathcal G$ on $X_2$. {\it
  ie.} there is a vector bundle $\mathcal G$ on $X_2$ such that the
multiplication map $f: \mathcal G(-d) \to \mathcal G$ induces the
exact sequence $0 \to G(-d) \to \mathcal G \to G \to 0$.
\end{lemma}

\begin{proof}
Since $Y$ is ACM, $H^1(I_{Y/X}(-d+\nu))=0, \forall \nu$, hence in the
sequence stated in the previous lemma, the right hand map is
surjective on the level of sections.  Therefore, the map
$\oplus_{i=1}^k\Oh_X(-m_i) \to I_{Y/X}$ can be lifted to a map
$\oplus_{i=1}^k\Oh_{X_2}(-m_i) \to I_{Y_2/X_2}$. Since a global
section of $I_{Y_2/X_2}(\nu)$ maps to zero in $I_{Y/X}$ only if it is
a multiple of $f$, by Nakayama's lemma, this lift is surjective at the
level of global sections in different twists, and hence on the level
of sheaves.  Hence there is a commuting diagram of exact sequences:

\[
\begin{array}{ccccccccc}
   0                &   &           0      & &   0        & & \\
 \uparrow           &   &       \uparrow  &  & \uparrow & & \\
 I_{Y_2/X_2}(-d)& \to &  I_{Y_2/X_2} & \to & I_{Y/X} &\to & 0 \\
\uparrow           &   &      \uparrow  & &   \uparrow & & \\
 
 \oplus_{i=1}^k\Oh_{X_2}(-m_i-d) & \to & \oplus_{i=1}^k\Oh_{X_2}(-m_i)
 & \to & \oplus_{i=1}^k\Oh_X(-m_i) & \to & 0\\

 \uparrow           &   &       \uparrow  & &  \uparrow & & \\ 
\mathcal G(-d)  &      \to &     \mathcal G &   \to &    G & \to & 0\\
 \uparrow           &   &       \uparrow  & &  \uparrow & & \\
   0                &   &           0    &  &   0        & & \\
\end{array}
\]

where the sheaf $\mathcal G$ is defined as the kernel of the lift, and
the map from the left column to the middle column is multiplication by
$f$. It is easy to verify that the lowest row induces an exact
sequence
$$ 0 \to G(-d) \to \mathcal G \to G \to 0.$$
By Nakayama's lemma,  $\mathcal G$ is a vector bundle on $X_2$. 
 
\end{proof}

\begin{prop}\label{split bundle}
Let $E$ be an ACM bundle on $X$. If $E$ extends to a bundle
$\mathcal{E}$ on $X_2$, then $E$ is a sum of line bundles.
\end{prop}

\begin{proof}
There is an exact sequence $0\to E(-d) \to \mathcal{E} \to E \to 0$,
where the left hand map is induced by multiplication by $f$ on
$\mathcal E$. Let $F_0=\oplus\Oh_{\bbP^n}(a_i) \onto E$ be a
surjection induced by the minimal generators of $E$. Since $E$ is ACM,
this lifts to a map $F_0\onto \mathcal{E}$. This lift is surjective on
global sections by Nakayama's lemma (since the sections of
$\mathcal{E}$ which are sent to $0$ in $E$ are multiples of $f$). Thus
we have a diagram
\[
\begin{array}{ccccccccc}
& &  & &  & & 0 & & \\
& &  & &  & & \downarrow & & \\
 & & 0 &  &  &  & E(-d)&  & \\
 & & \downarrow & &  & & \downarrow & & \\
 0 & \to & F_1 & \to & F_0 & \to & \mathcal{E} & \to & 0 \\
 & & \downarrow & & || & & \downarrow & & \\
0 & \to & G_1 & \to & F_0 & \to & E & \to & 0 \\
 & & \downarrow & &  & & \downarrow & & \\
  &  & E(-d) &  & &   & 0 \\
 & & \downarrow & &  & &  & & \\
 & & 0 & &  & &  & & \\
\end{array}
\]
$G_1$ and $F_1$ are sums of line bundles on $\bbP^n$ by Horrocks'
Theorem. Furthermore, $G_1\isom F_0(-d)$.  Thus $0 \to
F_0(-d)\by{\Phi} F_0 \to E \to 0$ is a minimal resolution for $E$ on
$\bbP^n$. As a consquence of this, one checks that
$\det{\Phi}=f^{\rank{E}}$. On the other hand, the degree of
$\det{\Phi}=d\rank{F_0}$ and so we have
$\rank{F_0}=\rank{E}$. Restricting, this resolution to $X$, we get a
surjection $F_0\tensor\Oh_X \onto E$. The ranks of both vector bundles
being the same, this implies that this is an isomorphism.
\end{proof}

\begin{cor}
Let $Y\subset X$ be a codimension $2$ ACM subvariety. If the conormal
sheaf sequence (\ref{nbs}) splits, then
\begin{itemize}
\item the ACM bundle $G$ associated to $Y$ is a sum of line bundles,
\item there is a codimension $2$ subvariety $S$ in $\bbP^n$ such that
$Y = X\cap S$.
\end{itemize}
\end{cor}

\begin{proof}

The first statement follows from Lemma \ref{bundle extends} and Proposition \ref{split bundle}.
For the second statement, since the bundle $G$ associated to $Y$ is a sum of line bundles 
$\oplus_{i=1}^{k-1}\Oh_X(-l_i)$ on $X$, the
map $G \to \oplus_{i=1}^k\Oh_X(-m_i)$ can be lifted to a map 
$\oplus_{i=1}^{k-1}\Oh_\bbP(-l_i) \to  \oplus_{i=1}^k\Oh_\bbP(-m_i)$. 
The determinantal variety $S$ of codimension $2$ in $\bbP^n$ determined 
by this map has the property that $Y = X\cap S$.
\end{proof}

In conclusion, we obtain the following collection of counterexamples:

\begin{cor}
If $G$ is an ACM bundle on $X$ which is not a sum of line bundles, and
if $Y$ is a subvariety of codimension $2$ in $X$ constructed from $G$
as in Proposition \ref{Kleiman}, then $Y$ does not satisfy the
conclusion of either Question \ref{ghconj1} or Question \ref{ghconj2}.
\end{cor}

Buchweitz-Greuel-Schreyer have shown \cite{BGS} that any hypersurface of 
degree at least $2$ supports (usually many) non-split ACM bundles. We will 
give another construction in the next section.

\section {Remarks}

\begin{subsection}{}
{The infinitesimal Question \ref{ghconj2} was treated by studying the
  extension of the bundle to the thickened hypersurface $X_2$. This
  method goes back to Ellingsrud, Gruson, Peskine and Str{\o}mme
  \cite{EGPS}.}  If we are not interested in the infinitesimal
Question \ref{ghconj2}, but just in the more geometric Question
\ref{ghconj1}, a geometric argument gives an even easier proof of the
existence of codimension $2$ ACM subvarieties $Y\subset X$ which are
not of the form $Y=X\cap Z$ for some codimension $2$ subvariety
$Z\subset \bbP^n$.

\begin{prop}
Let $E$ be an ACM bundle on a hypersurface $X$ in $\bbP^n$ which
extends to a sheaf $\mathcal{E}$ on
\comment{\textcolor[rgb]{0.98,0.00,0.00}}$\bbP^n$; i.e. there is an exact
sequence
\begin{equation}\label{extension}
0 \to \mathcal{E}(-d) \stackrel{f}{\to} \mathcal{E} \to E \to 0  
\end{equation}

Then $E$ is a sum of line bundles.
\end{prop}

\begin{proof}
At each point $p$ on $X$, over the local ring $\Oh_{\bbP,p}$ the sheaf
$\mathcal E$ is free, of the same rank as $E$. Hence $\mathcal E$ is
locally free except at finitely many points. Let $\bbH$ be a general
hyperplane not passing through these points. Let $X' = X\cap \bbH$,
and $\mathcal E', E'$ be the restrictions of $\mathcal E, E$ to $\bbH,
X'$.

It is enough to show that $E'$ is a sum of line bundles on $X'$. This
is because any isomorphism $\oplus \Oh_{X'}(a_i) \to E'$ can be lifted
to an isomorphism $\oplus\Oh_{X}(a_i) \to E$, as $H^1(E(\nu))=0,
\forall~ \nu\in\mathbb{Z}$. The bundle $E'$ on $X'$ is ACM and from
the sequence
$$ 0 \to \mathcal E'(-d) \to \mathcal E' \to E' \to 0,$$ it is easy to
check that $H^i(\mathcal E'(\nu))=0, \forall~ \nu\in\mathbb{Z}$, for
$2 \leq i \leq n-2$. Since $\mathcal E'$ is a vector bundle on $\bbH$,
we can dualize the sequence to get
$$ 0 \to \mathcal E'^{\vee}(-d) \to \mathcal E'^{\vee} \to E'^{\vee}
\to 0.$$ $E'^{\vee}$ is still an ACM bundle, hence $H^i(\mathcal
E'^{\vee}(\nu)) =0, \forall~ \nu\in\mathbb{Z}$, and $2 \leq i \leq
n-2$.

By Serre duality, we conclude that $\mathcal E'$ is an ACM bundle on $\bbH$,
and by Horrocks' theorem, $\mathcal E'$ is a sum of line bundles. Hence, its
restriction $E'$ is also a sum of line bundles on $X'$.
\end{proof}

\begin{prop}
Let $Y$ be an ACM subvariety of codimension $2$ in the hypersurface $X$ such 
that the associated ACM bundle $G$ is not a sum of line bundles.
Then there is no pure subvariety $Z$ of codimension $2$ in $\bbP^n$ such that
$Z\cap{X}=Y$. 
\end{prop}

\begin{proof}
Suppose there is such a $Z$. Then there is an exact sequence $0 \to
{I}_{Z/\bbP}(-d) \to {I}_{Z/\bbP} \to {I}_{Y/X} \to 0$, where the
inclusion is multiplication by $f$, the polynomial defining $X$. Since
$Z$ has no embedded points, $H^1({I}_{Z/\bbP}(\nu))=0$ for $\nu<<0$.
Combining this with $H^1({I}_{Y/X}(\nu))=0, \forall~\nu\in\bbZ$, and using
the long exact sequence of cohomology, we get
$H^1({I}_{Z/\bbP}(\nu))=0, \forall~ \nu\in\bbZ$.

Now suppose $Y$ has the resolution $ 0\to G \to \oplus\Oh_{X}(-m_i)\to
I_{Y/X} \to 0$. From the vanishing just proved, the right hand map can
be lifted to a map $\oplus\Oh_{\bbP}(-m_i)\to I_{Z/\bbP}$, which is
easily checked to be surjective (at the level of global sections). It
follows that if $\mathcal G$ is the kernel of this lift, $\mathcal G$
is an extension of $G$ to $\bbP^n$. By the previous proposition, $G$
is a sum of line bundles. This is a contradiction.
\end{proof}

\end{subsection}

\begin{subsection}{}
Voisin's original example was as follows. Let $P_1$ and $P_2$ be two
planes meeting at a point $p$ in $\bbP^4$. The union $\Sigma$ is a
surface which is not locally Cohen-Macaulay at $p$. Let $X$ be a
smooth hypersurface of degree $d>1$ which passes through $p$. $X\cap
\Sigma$ is a curve $Z$ in $X$ with an embedded point at $p$. The
reduced subscheme $Y$ has the form $Y=C_1\cup C_2$, where $C_1$ and
$C_2$ are plane curves. Voisin argues that $Y$ itself does not have
the form $X\cap S$ for any surface $S$ in $\bbP^4$.

We can treat this example from the point of view of ACM bundles.
$I_{Z/X}$ has a resolution on $X$ which is just the restriction of the
resolution of the ideal of the union $P_1\cup P_2$ in $\bbP^4$, {\it
  viz.}
$$ 0 \to \Oh_X(-4) \to 4\Oh_X(-3) \to 4 \Oh_X(-2) \to I_{Z/X} \to
0. $$ From the sequence $0 \to I_{Z/X} \to I_{Y/X} \to k_p \to 0$, it
is easy to see that $Y$ is ACM, with a resolution
$$ 0 \to G \to 4\Oh_X(-2)\oplus \Oh_X(-d) \to I_{Y/X} \to 0.$$ $G$ is
an ACM bundle. If it were a sum of line bundles, comparing the two
resolutions, we find that $h^0(G(2)) =0$ and $h^0(G(3)) =4$, hence $G=
4\Oh_X(-3)$. But then $G \to 4\Oh_X(-2)\oplus \Oh_X(-d)$ cannot be an
inclusion. Thus $G$ is an ACM bundle which is not a sum of line
bundles.

Voisin's subsequent smooth examples were obtained by placing $Y$ on a
smooth surface $T$ contained in $X$ and choosing divisors $Y'$ in the
linear series $|Y+mH|$ on $T$. When $m$ is large, $Y'$ can be chosen
smooth. In fact, such curves $Y'$ are doubly linked to the original
curve $Y$ in $X$, hence they have a similar resolution $ G' \to L \to
I_{D'/X}\to 0$, where $L$ is a sum of line bundles and where $ G'$
equals $G$ up to a twist and a sum of line bundles.

The fact that $G$ above is not a sum of line bundles is related 
(via the mapping cone of the map of resolutions)
to  the fact that $k_p$ itself cannot have a finite resolution by sums 
of line bundles on $X$. This follows from the
following proposition which provides another argument for the existence of ACM 
bundles on arbitrary smooth hypersurfaces of degree $\geq 2$.

\begin{prop}\label{constructions}
Let $X$ be a smooth hypersurface in $\bbP^n$ of degree $\geq 2$ with
homogeneous coordinated ring $S_X$. Let $L$ be a linear space
(possibly a point or even empty) inside $X$ of codimension $r$, with
homogeneous ideal $I(L)$ in $S_X$.  A free presentation of $I(L)$ of
length $r-2$ will have a kernel whose sheafification is an ACM bundle
on $X$ which is not a sum of line bundles.
\end{prop}

\begin{proof} It should first be understood that the homogeneous ideal $I(L)$ 
of the empty linear space will be taken as the irrelevant ideal 
$(X_0, X_1, \dots, X_n)$. Let the free presentation of $I(L)$ together with the
kernel be
$$ 0\to M \to F_{r-2} \to \cdots \to F_0 \to I(L) \to 0, $$ 
where $F_i$ are free graded $S_X$ modules. Its sheafification looks like
$$ 0 \to \tilde M \to \tilde F_{r-2} \to \cdots \to \tilde F_0 \to
I_{L/X} \to 0.$$ Since $L$ is locally Cohen-Macaulay, $\tilde M$ is a
vector bundle on $X$, and since $L$ is ACM, so is $\tilde M$. $M$
equals $\oplus_{\nu \in \bbZ} H^0(\tilde M(\nu))$.  Hence, $\tilde M$
is a sum of line bundles only if $M$ is a free $S_X$ module.

If $\bbH$ is a general hyperplane in $\bbP^n$ which meets $X$ and $L$
transversally along $X_{\bbH}$ and $L_{\bbH}$ respectively, the above
sequences of modules and sheaves can be restricted to give similar
sequences in $\bbH$. The restriction $\tilde M_{\bbH}$ is an ACM
bundle on $X_{\bbH}$.

Repeat this \comment{\textcolor[rgb]{0.98,0.00,0.00}}successively to
find a maximal and general linear space $\bbP$ in $\bbP^n$ which does
not meet $L$.  If $X' = X\cap \bbP$, the restriction of the sequence
of $S_X$ modules to $X'$ gives a resolution
$$ 0 \to M' \to F_{r-2}' \to \cdots \to F_0' \to S_{X'} \to k \to 0.$$

Localize this sequence of graded $S_{X'}$ modules at the irrelevant
ideal $I(L)\cdot S_{X'}$, to look at its behaviour at the vertex of
the affine cone over $X'$. $k$ is the residue field of this local
ring. Since $X$ and hence $X'$ has degree $\geq 2$, the cone is not
smooth at the vertex. By Serre's theorem (\cite{Se}, IV-C-3-Cor 2),
$k$ cannot have finite projective dimension over this local
ring. Hence $M'$ is not a free module. Therefore neither is $M$.
\end{proof}
\end{subsection}

\begin{subsection}{}
We make a few concluding remarks about Question \ref{ghconj3}, the
Degree Conjecture of Griffiths and Harris. A vector bundle $G$ on a
smooth hypersurface $X$ in $\bbP^4$ has a second Chern class $c_2(G)
\in A^2(X)$, the Chow group of codimension $2$ cycles.  If $h \in
A^1(X)$ is the class of the hyperplane section of $X$, the degree of
any element $c\in A^2(X)$ will be defined to be the degree of the zero
cycle $c\cdot h \in A^3(X)$. (Note that by the Lefschetz theorem, all
classes in $A^1(X)$ are multiples of $h$.)

With this notation, if $E$ is any bundle on $X$ and $Y$ is a curve
obtained from $E$ with the sequence ({\it vide} Proposition
\ref{Kleiman})
$$ 0 \to E \to \oplus_{i=1}^k\Oh_X(m_i) \to I_{Y/X}(m) \to 0,$$ a
calculation tells us that the degree $d$ of $X$ divides the degree of
$Y$ if and only if $d$ divides the degree of $c_2(E)$.

More generally: let $Y$ be any curve in $X$ and resolve $I_{Y/X}$ to get
$$ 0 \to E \to \oplus_{i=1}^l\Oh_X(b_i) \to \oplus_{i=1}^k\Oh_X(a_i)
\to I_{Y/X} \to 0,$$

where $E$ is an ACM bundle on $X$. Then a similar calculation tells us
that the degree $d$ of $X$ divides the degree of $Y$ if and only if
$d$ divides the degree of $c_2(E)$.

Hence we may ask the following question which is equivalent to the
Degree Conjecture:

\begin{ACM} If $X$ is a general hypersurface in $\bbP^4$ of degree
 $d \geq 6$, then for any indecomposable ACM vector bundle $E$ on $X$,
  $d$ divides the degree of $c_2(E)$.
\end{ACM}

The examples created above in Proposition \ref{constructions} satisfy
this, when $L$ has codimension $> 2$ in $X$.
In \cite{MRR}, this conjecture is settled for ACM bundles of rank $2$ on $X$.

\end{subsection}

\bibliographystyle{alpha}

\end{document}